\documentclass[12pt]{amsart}
\usepackage{amsmath, amssymb, latexsym}

\usepackage[dvips]{graphicx}
\usepackage[dvips]{color}

\title{On the Cauchy transform of the Bergman space}
\author{S. A. Merenkov}
\thanks{Research supported in part by INTAS, reference number INTAS 96-0858}

\newtheorem{theorem}{Theorem}
\newtheorem{definition}{Definition}
\newtheorem{lemma}{Lemma}
\newtheorem{remark}{Remark}

\begin{document}

\abstract{The range of the Bergman space $B_2(G)$ under the Cauchy transform $K$
is described for a large class of domains. For a quasidisk $G$ the
relation $K(B_2^*(G))=B_2^1(\Bbb C\setminus\overline{G})$ is proved.}

\endabstract

\maketitle

\section{Introduction}

Let $G$ be a domain in the complex plane $\Bbb C$ bounded by a Jordan curve $\partial\, G$ with area($\partial\, G$)=0. We call these domains integrable domains. Consider
the following classes of analytic functions:
$$
B_2(G)=\bigg\{ g(z)\in {\text {Hol}}(G),{\| g\| }_{B_2(G)}=\biggl(\iint\limits_G|g(z)|^2dxdy
\biggr)^{\frac12}<\infty \bigg\} ;
$$
$$
H(\Bbb C\setminus\overline G)=\{\gamma(\zeta)\in {\text {Hol}}(\Bbb C\setminus\overline G),
\gamma(\infty)=0\};
$$
$$
\begin{aligned}
B_2^1(\pmb{\Bbb C}\setminus\overline G)=\bigg\{&\gamma(\zeta)\in H(\Bbb C\setminus\overline G),\\
&{\|\gamma \| }_{B_2^1(\pmb{\Bbb C}\setminus\overline G)}=\biggl(\iint\limits_
{\pmb{\Bbb C}\setminus \overline G}|\gamma'(\zeta)|^2d\xi d\eta \biggr)^{\frac12
}<\infty\bigg\},
\end{aligned}
$$
where $z=x+iy,\,\zeta=\xi+i\eta;$
$\overline G$ is the
closure of the domain $G$. The class $B_2(G)$ is called the Bergman space.
 
The transformation
$$
(Kg)(\zeta)=\frac 1\pi \iint \limits_G \frac{ \overline {g(z)}}{z-\zeta}\,dxdy,
$$
where $g(z)\in B_2(G),\,\zeta \notin \overline G$ is called the Cauchy 
transform of $B_2^\ast(G)$ which is dual to $B_2(G)$. Because the spaces $B_2(G)$ and 
$B_2^\ast(G)$
are isometric, we can think of $K$ as a transformation of $B_2(G)$.

The problem of describing the range of $X^{\ast}$ under the 
Cauchy transform for different spaces $X$ of analytic functions was investigated by many authors,
see, for example, [6]-[7]. The motivation 
of the present work is the paper [1]. V.V.Napalkov(jr) and R.S.Yul\-mu\-kha\-me\-tov established that $K(B_2^{\ast}(G))=B_2^1(\Bbb C\setminus\overline G)$ for domains
with sufficiently smooth boundary. We prove that this relation is valid
for quasidisks, and also find $K(B_2^{\ast}(G))$ for a large class
of domains.


It is obvious that the Cauchy transform converts a function $g(z)\in B_2(G)$ 
into an analytic function $\gamma(\zeta)$ on
$\Bbb C\setminus{\overline G}$ such that $\gamma(\infty)=0$. Since polynomials
are dense in $B_2(G)$ [2, ch.1, 3] and the system $\{1/(z-\zeta),\,\, \zeta\notin\overline G\}$ is dense in the space of functions holomorphic in $\overline G$, the operator $K$ is injective.
%

The operator
$$
(\Bbb Tu)(\zeta)=\frac1{\pi}\lim_{\varepsilon \to 0}
\iint\limits_{|z-\zeta|\geq\varepsilon}\frac{u(z)}{(z-\zeta)^2}\,dxdy
$$
is an isometry on $L_2(\Bbb C)$ [3, pp. 64-66]. Thus 
$K:B_2^{\ast}(G)\rightarrow B_2^1(\Bbb C\setminus\overline G)$
is a continuous operator.


Throughout the paper we denote the unit disk by $\Bbb D$ and its boundary by $\partial\, \Bbb D$. The boundary of a domain $G$ is denoted by $\partial\, G$.

{\bf{Acknowledgment:}} 
The author gratefully thanks A.F.Grishin and M.Sodin for their interest and valuable discussions.

\section{General case}

To study $K(B_2^{\ast}(G))$ we need the function space
$$
\gathered
W(0,2\pi)=\left\{f(e^{i\theta})\in L_1(0,2\pi),\,f(e^{i\theta})\sim
\sum_{k=-\infty}^{\infty}{f_ke^{ik\theta}},\right. \\
\left. {\text {with the semi-norm}}\quad \rho(f)=\biggl(\pi\sum_{k=1}^{\infty}{k|f_{-k}|^2}\biggr)^{\frac12}<\infty \right\}.
\endgathered
$$

Functions of $W(0,2\pi)$ can be characterized as follows:

\begin{lemma}
Let $f(t)\in L_1(\partial\,\Bbb D),$ i.e. $f(e^{i\theta})\in L_1(0,2\pi),$
and $F(\zeta)$ be the Cauchy-type integral corresponding to $f(t)$:
$$
F(\zeta)=\frac1{2\pi i}\int\limits_{\partial\,\Bbb D}\frac{f(t)}{t-\zeta}\,dt
,\quad\zeta \in\Bbb C\setminus{\overline{\Bbb D}}.
$$
Then $f\in W(0,2\pi)$ if and only if $F 
\in B_2^1(\Bbb C\setminus\overline{\Bbb D}),$ and
$$
\rho(f)=\|F\|_{B_2^1(\Bbb C\setminus\overline{\Bbb D})}.
$$ 
\end{lemma}

\emph{Proof}
It is obvious that $F(\zeta)\in {\text {Hol}}(\Bbb C\setminus\overline G)$ and
$F(\infty)=0.$

Next we have 
%
$$
\gathered
F(\zeta)=\frac1{2\pi i}\int\limits_{\partial\, \Bbb D}\frac{f(t)}{t-\zeta}\,dt=-\frac1{\zeta}
\sum\limits_{k=0}^{\infty}\frac1{\zeta^k}\frac1{2\pi i}\int\limits_{\partial\, \Bbb D}
f(t)\,t^kdt=-\sum\limits_{k=1}^{\infty}\frac{f_{-k}}{\zeta^k}.
\endgathered
$$
The identity $\quad\|F\|_{B_2^1(\Bbb C\setminus\overline{\Bbb D})}=\bigg(\pi\sum\limits_{k=1}^{\infty}|F_k|^2k\bigg)^{\frac12},
$ where $\{F_k\}_1^{\infty}$ is the set of Taylor coefficients of $F$, proves the lemma.

\qed


Let $G$ be an integrable domain and let a sequence of Jordan domains  $\{G_n\}_1^{\infty}$ 
satisfies the conditions:
 
(i)  $\partial G_n$
 is a smooth Jordan curve; (ii) $\overline G_{n+1}\subset G_n,n=1,2,3,\dots$;
(iii) $ \cap_{n\geq 1}G_n=\overline G$.
Let $\varphi_n$ be a conformal map of $\Bbb D$ onto $G_n$.
\begin{theorem}
A function $\gamma$ from $ B_2^1(\Bbb C\setminus\overline G)$ belongs to $K(B_2^{\ast}(G))$
if and only if 
$$
\sup_{n\geq1}\rho(\gamma\circ\varphi_n(e^{i\theta}))<\infty
$$
for any sequence $\{G_n\}_1^{\infty}$ with (i),(ii),(iii).
\end{theorem}
\emph{Proof}
First we show the implication
$$
\left\{\gamma(\zeta)\in B_2^1(\Bbb C\setminus\overline G),\quad
\sup_{n\geq1}\rho(\gamma\circ\varphi_n(e^{i\theta}))<\infty\right\}\subset
K(B_2(G)).
$$

Let $\gamma(\zeta)$ belong to $B_2^1(\Bbb C\setminus\overline G)$ and
$\sup_{n\geq1}\rho(\gamma\circ\varphi_n(e^{i\theta}))<\infty$.
We write $h\in {\text {Hol}}(\overline G)$ if there exists an open set $G_1=G_1(h)\supset
\overline G$ such that $h\in {\text {Hol}}(G_1)$. For functions $h\in {\text {Hol}}(\overline G)$
 introduce the linear functional:
$$
\Bbb F(h)=\lim_{n\to\infty}\int\limits_{\partial G_n}\gamma (\xi)\,h(\xi)\,d\xi.
$$
If $n_0$ is  such a number that $h$ is holomorphic in $G_{n_0}$, then the last
integral is unaffected by $n\geq n_0$. Thus, $F(h)$ is meaningful.

We show that $\Bbb F$ is a bounded linear functional
on the space ${\text {Hol}}(\overline G)$ using the norm of the space
$B_2(G)$. Changing the variable by formula $\xi=\varphi_n(e^{i\theta})$,
 get
$$
\frac1{2\pi i}\int\limits_{\partial G_n}\gamma(\xi)h(\xi)\,d\xi=
\frac1{2\pi i}\int\limits_0^{2\pi}\gamma(\varphi_n(e^{i\theta}))
h(\varphi_n(e^{i\theta}))(\varphi_n)'_{\theta}(e^{i\theta})\,d\theta.
$$
The function $h(\varphi_n(e^{i\theta}))(\varphi_n)'_{\theta}(e^{i\theta})$
is the restriction to the unit circumference of the function
$h_n(z)=h(\varphi_n(z))(\varphi_n)'(z)zi$
[4, p.405].
%
%
Changing the variable $w=\varphi_n(z)$ we see that $\|h_n\|_{B_2(\Bbb D)}\leq\|h\|_{B_2(G_n)}$.
Since $h(z)$ is continuous in $\overline G_n$ for $n\geq n_0$ and $\varphi_n(z)$
maps the unit disk onto the domain $G_n$ bounded by a smooth
Jordan curve, $\varphi_n'(z)$ and $h_n(z)$ belong to $H_2(\Bbb D)$
(Hardy space) [4, p.410]. If $\{c_k^n\}_1^{\infty}$ is the sequence
of Taylor coefficients for the function $h_n(z)$, then an easy calculation
shows


$$
\|h_n\|_{B_2(\Bbb D)}=\bigg(\pi\sum\limits_{k=1}^{\infty}
\frac{|c_k^n|^2}{k+1}\bigg)^{\frac12}<\infty.
$$
Thus
$$\frac1{2\pi i}\int\limits_{\partial G_n}\gamma(\xi)h(\xi)\,d\xi=
\frac1{2\pi i}\int\limits_0^{2\pi}\gamma(\varphi_n(e^{i\theta}))h_n(e^{i\theta})
\,d\theta=\frac1i\sum\limits_{k=1}^{\infty}a_{-k}^nc_k^n,
$$
where
$\{a_n^k\}_{-\infty}^{\infty}$ is defined by the formula
$\gamma(\varphi_n(e^{i\theta}))=\sum\limits_{k=-\infty}^{\infty}a_k^ne^{ik\theta}$.
Applying the Cauchy-Schwarz inequality, we get
$$
\gathered
|\frac1{2\pi i}\int\limits_{\partial G_n}\gamma(\xi)h(\xi)\,d\xi|=
|\sum\limits_{k=1}^{\infty}a_{-k}^nc_k^n|
\leq\bigg(\sum\limits_{k=1}^{\infty}k|a_{-k}^n|^2\bigg)^{\frac12}
\bigg(\sum\limits_{k=1}^{\infty}\frac{|c_k^n|^2}{k}\bigg)^{\frac12}\\
=\frac{\sqrt2}{\pi}\rho(\gamma\circ\varphi_n(e^{i\theta}))\|h_n\|_{B_2(\Bbb D)}
\leq \frac{\sqrt2}{\pi}\rho(\gamma\circ\varphi_n(e^{i\theta}))\|h\|_{B_2(G_n)}.
\endgathered
$$
%
Because the domain $G$ is
integrable, $\|h\|_{B_2(G_n)}\rightarrow\|h\|_{B_2(G)}$ as
$n\rightarrow\infty$. Hence
$$
|\Bbb F(h)|\leq C\|h\|_{B_2(G)},\qquad \text {where}\quad
C=\frac{\sqrt2}{\pi}\sup\limits_{n\geq1}\rho(\gamma\circ\varphi_n(e^{i\theta})).
$$
Since the space $\text {Hol}(\overline G)$ is dense in $B_2(G)$, the functional $\Bbb F$
can be uniquely extended to the linear continuous functional
on $B_2(G)$ that we denote by $F$ also. It follows from the Riesz-Fisher representation theorem
that there exists a function $g\in B_2(G)$ such that
$$
\Bbb F(h)=\frac 1{\pi}\iint\limits_Gh(z)\overline{g(z)}\,dxdy,\quad
h\in B_2(G).
$$
Now calculate $\Bbb F(1/(z-\zeta))$
for $\zeta\notin\overline G$,
$$
\Bbb F\big(\frac1{z-\zeta}\big)=\lim\limits_{n\rightarrow\infty}\frac1{2\pi i}\int\limits_{\partial G_n}
\frac{\gamma(z)}{z-\zeta}\,dz=-\gamma(\zeta).
$$
We obtain that
$$
\gamma(\zeta)=\frac1{\pi}\iint\limits_G\frac{-\overline{g(z)}}{z-\zeta}\,dxdy,
\quad \zeta\notin\overline G\quad \text{ and }-g\in B_2(G).
$$
The relation 
$$
\left\{\gamma(\zeta)\in B_2^1(\Bbb C\setminus\overline G),\quad
\sup_{n\geq1}\rho(\gamma\circ\varphi_n(e^{i\theta}))<\infty\right\}\subset
K(B_2(G))
$$
is proved.

To prove the relation
$$
K(B_2(G))\subset\left\{\gamma(\zeta)\in B_2^1(\Bbb C\setminus\overline G),\quad
\sup_{n\geq1}\rho(\gamma\circ\varphi_n(e^{i\theta}))<\infty\right\}
$$
we  apply the lemma. It is sufficient to show that
$\sup_{n\geq 1}\|F_n\|_{B_2^1(\Bbb C\setminus \overline{\Bbb D})}<\infty$,
where
$$
F_n(\zeta)=\frac1{2\pi i}\int\limits_{\partial\, \Bbb D}\frac{\gamma\circ\varphi_n
(t)}{t-\zeta}\,d\zeta,\quad \gamma(\zeta)=\frac1\pi\iint\limits_G\frac{\overline{g(z)}}{z-\zeta}\,dxdy,
\quad g\in B_2(G).
$$
Putting the expression for $\gamma(\zeta)$ in the formula for $F_n(\zeta)$,
 have
$$
F_n(\zeta)=\frac 1{2\pi i}\int\limits_{\partial\, \Bbb D}\frac1{t-\zeta}\frac1\pi\iint\limits_G
\frac{\overline{g(z)}}{z-\varphi_n(t)}\,dxdy\,dt.
$$
Since $\overline{g(z)}/((t-\zeta)(z-\varphi_n(t)))\in L_1(G\times
\partial\, \Bbb D)$ for $\zeta\in\Bbb C\setminus\overline{\Bbb D}$,  we can
 interchange the order of integration: 
$$
F_n(\zeta)=\frac1\pi\iint\limits_G\overline{g(z)}\frac1{2\pi i}\int \limits_{\partial\,
\Bbb D}\frac1{t-\zeta}\,\frac1{z-\varphi_n(t)}\,dt\,dxdy.
$$
Further, the residue theorem yields
$$
\frac1{2\pi i}\int\limits_{\partial\,\Bbb D}\frac1{t-\zeta}\,\frac1{z-\varphi_n(t)}\,dt
=-\frac 1{(\varphi_n^{-1}(z)-\zeta)\varphi_n'(\varphi_n^{-1}(z))},
$$
where $\varphi_n^{-1}$ is the inverse function of $\varphi_n$.
Let $w=\varphi_n^{-1}(z)$ in the resulting integral, we then see that 
$$
F_n(\zeta)=-\frac 1{\pi} \iint\limits_{\Bbb D_n}\frac {\overline{g(\varphi_n(w))
\varphi_n'(w)}}{w-\zeta}\,dudv,
$$
where $\Bbb D_n=\varphi_n^{-1}(G) \subset\Bbb D$.
Hence in $\Bbb C\setminus \overline{\Bbb D}$
$$
F_n'(\zeta)=\Bbb T(-\overline{g
(\varphi_n(w))\varphi_n'(w)})(\zeta),
$$
where the operator $\Bbb T$  was introduced earlier.
Since $\Bbb T$ is isometric, we get
$$
\gathered
\|F_n\|_{B_2^1(\Bbb C\setminus\overline{\Bbb D})}\leq\|\Bbb T(-\overline{g
(\varphi_n(w))\varphi_n'(w)})\|_{L_2(\Bbb C\setminus\overline{\Bbb D})}\\
\leq \|g(\varphi_n(w))\varphi_n'(w)\|_{B_2(\Bbb D_n)}=\|g\|_{B_2(G)}.
\endgathered
$$
Thus
$$
\sup\limits_{n\geq1}\|F_n\|_{B_2^1(\Bbb C\setminus\overline{\Bbb D})}\leq \|g\|_{B_2(G)},
$$
Theorem 1 is proved.
\qed
\section{The case of a quasidisk}

As an application of Theorem 1  we prove a theorem concerning the 
Cauchy transform of the Bergman space on quasidisks.

We give some definitions [5, ch. 5].

\begin{definition} A quasiconformal map of $\Bbb C$ onto $\Bbb C$
is a homeomorphism $h$ such that:
\begin{enumerate}
\item $h(x+iy)$ is absolutely continuous in $x$ for almost all $y$ and in 
$y$ for almost all $x$;
\item the partial derivatives are locally square integrable;
\item $h(x+iy)$ satisfies the Beltrami differential equation
\end{enumerate}
$$
\frac{\partial h}{\partial\overline z}=\mu(z)\frac{\partial h}{\partial z}
\quad\text{for almost all }z\in\Bbb C,
$$
where  $\mu$ is a complex measurable function with $|\mu(z)|\leq k<1$
for $z\in\Bbb C$.
In this case it is said  $h$ to be a $k-$quasiconformal map.
\end{definition}

\begin{definition} A quasicircle in $\Bbb C$ is a Jordan curve $J$ such that
$$
{\text {diam}} J(a,b)\leq M|a-b|\quad \text {for}\quad a,b\in J,
$$
where $J(a,b)$ is the arc of the smaller diameter  of $J$ between $a$ and $b$. The 
domain interior to $J$ is called a quasidisk.
\end{definition}
\begin{remark}
An equivalent definition for $J$ to be a quasicircle: $J$ is the range 
of the circle under a quasiconformal map of $\Bbb C$ onto $\Bbb C$.
\end{remark}

\begin{theorem}
Let $G$ be a quasidisk, then
$$
K(B_2^{\ast}(G))=B_2^1(\Bbb C\setminus\overline G).
$$
\end{theorem}
\emph{Proof}
Let $\psi$ be a conformal map of $\Bbb C\setminus\overline{\Bbb D}$
onto $\Bbb C\setminus\overline G$ with $\psi(\infty)=\infty$. Denote the inner
domain bounded by the curve $\{\psi(R_ne^{i\theta}), \theta\in[0,2\pi)\}$ by
$G_n$, where $\{R_n\}_1^{\infty}$ be some sequence decreasing monotonically to 1.
Let $\varphi_n$ be a conformal map of $\Bbb D$ onto $G_n$.

Since $K(B_2^{\ast}(G))\subset B_2^1(\Bbb C\setminus\overline G)$, we have only to show
that for every $\gamma(\zeta)\in B_2^1(\Bbb C\setminus\overline G)$
the following holds true: $\sup_{n\geq1}\rho(\gamma\circ\varphi_n(e^{i\theta}))
<\infty$.
Then, in view of Theorem 1, we get Theorem 2.

To verify the inequality $\sup_{n\geq1}\rho(\gamma\circ\varphi_n(e^{i\theta}))<\infty$
 apply the lemma. We have
$$
F_n(\zeta)=\frac1{2\pi i}
\int\limits_{\partial\,\Bbb D}\frac{\gamma\circ\varphi_n(t)}{t-\zeta}dt,\quad
\zeta\in\Bbb C\setminus\overline{\Bbb D}.
$$
It is clear that 
$$
|\psi^{-1}\circ\varphi_n(t)|=R_n,\quad t\in\partial\,\Bbb D,\quad n\geq 1.
$$
Hence $\gamma\circ\varphi_n(t)=
\gamma\circ\psi\left({R_n^2}/({\overline{\psi^{-1}\circ\varphi_n(t)}})\right),
t\in\partial\,\Bbb D$.

Theorem 5.17 [5, p.114] states that any conformal map of the disk
onto a quasidisk can be extended to a quasiconformal map of $\Bbb C$ onto $\Bbb C$.
Evidently, the theorem remains true for a conformal map of $\Bbb C\setminus
\overline{\Bbb D}$ onto a domain exterior to a quasicircle. It gives that the function $\psi$
can be extended to a quasiconformal map $\Psi:\quad\Bbb C\rightarrow\Bbb C$.
Let $\Psi$ be a $k-$quasiconformal map. Then $\Psi^{-1}$ is of that kind.
Composition of a conformal and a $k-$quasiconformal maps is $k-$quasiconformal. 
Thus the function $\overline{f_n}(z)={R_n^2}/(
{\Psi^{-1}\circ\varphi_n(z)})$ is $k-$quasiconformal map of  $\Bbb D$  onto
$\{|w|>R_n\}$, 
$\left|\partial f_n/\partial z\right|\leq k\left|\partial f_n/\partial\overline z\right|.$ 
If $J_n$ stands for the Jacobian of $f_n, J_n=\big|
\partial f_n/\partial z\big|^2-\big|\partial f_n/\partial\overline z
\big|^2$, then $\big|\partial f_n/\partial z\big|^2\leq|J_n|/(1-k^2).
$

We need to estimate $\big\|\partial/\partial\overline z\,\gamma\circ
\psi(f_n(z))\big\|_{L_2(\Bbb D)}$.
$$
\gathered
\big\|\frac\partial{\partial\overline z}\gamma\circ\psi(f_n(z))\big\|_{L_2(\Bbb D)}
=\bigg(\iint\limits_{\Bbb D}|(\gamma\circ\psi)'(f_n(z))|^2|\frac{\partial}
{\partial\overline z}f_n(z)|^2dxdy\bigg)^{\frac12}\\
\leq\frac1{\sqrt{1-k^2}}\bigg(\iint\limits_{\Bbb D}|(\gamma\circ\psi)'(f_n(z))|^2|J_n(z)|\,dxdy\bigg)^{\frac12}
\endgathered
$$
$$
 \leq\frac1{\sqrt{1-k^2}}\bigg(\iint\limits_{\Bbb C\setminus\overline{\Bbb D}}|(\gamma\circ\psi)'(w)|^2dudv\bigg)
^{\frac12},
$$
where $w=u+iv$.
Since the operator $\tilde \psi:\,\, \tilde \psi(\gamma)(\zeta)=\gamma\circ\psi(\zeta)$ 
is an isometry from $B_2^1(\Bbb C\setminus \overline G)$ to $B_2^1(\Bbb C\setminus
\overline{\Bbb D})$,
we have:
$$
\big\|\frac\partial{\partial\overline z}\gamma\circ\psi(f_n(z))\big\|_{L_2(\Bbb D)}
\leq\frac1{\sqrt{1-k^2}}\|\gamma\|_{B_2^1(\Bbb C\setminus\overline G)}.
$$

Now the Green formula gives:
$$
F_n(\zeta)=\frac1{2\pi i}\int\limits_{\partial\,\Bbb D}\frac{\gamma\circ\psi(f_n(t))}{t-\zeta}dt=
\frac1{\pi}\iint\limits_{\Bbb D}\frac1{z-\zeta}\frac\partial
{\partial\overline z}\gamma\circ\psi(f_n(z))\,dxdy
$$
for $\zeta\in\Bbb C\setminus\overline{\Bbb D}$.
Using isometricity of the operator $\Bbb T$ defined above, we get
$$
\|F_n\|_{B_2^1(\Bbb C\setminus\overline{\Bbb D})}\leq \big\|\frac
\partial{\partial\overline z}
\gamma \circ\psi(f_n(z))\big\|_{L_2(\Bbb D)}\leq
\frac1 {\sqrt{1-k^2}}\|\gamma\|_{B_2^1(\Bbb C\setminus\overline G)}.
$$
Thus Theorem 2 is proved.

\qed


\end{document}